\documentclass[a4paper,11pt,leqno]{amsart}

\usepackage{amsmath,amssymb,array,euscript,amsfonts}
\usepackage{graphicx,epsfig,color}
\usepackage[usenames,dvipsnames]{xcolor}
\usepackage{rotating}
\usepackage{color}

\def\C{\hbox{\font\dubl=msbm10 scaled 1100 {\dubl C}}}
\def\R{\hbox{\font\dubl=msbm10 scaled 1100 {\dubl R}}}

\def\N{\hbox{\font\dubl=msbm10 scaled 1100 {\dubl N}}}

\def\sR{\hbox{\font\dubl=msbm10 scaled 900 {\dubl R}}}

\def\d{\,{\rm{d}}}
\def\Re{{\rm{Re}}\,}
\def\Im{{\rm{Im}}\,}

\newtheorem{Theorem}{Theorem}

\title[Riemann Zeta-Function along its Julia Lines]{Value-Distribution of the Riemann Zeta-Function along its Julia Lines}

\thanks{{\bf Dedicated to the Memory of Stephan Ruscheweyh}}

\author[J\"orn Steuding and Ade Irma Suriajaya]{J\"orn Steuding and Ade Irma Suriajaya}

\date{November 2019}

\begin{document}

\maketitle

\begin{abstract} 
For an arbitrary complex number $a\neq 0$ we consider the distribution of values of the Riemann zeta-function $\zeta$ at the $a$-points of the function $\Delta$ which appears in the functional equation $\zeta(s)=\Delta(s)\zeta(1-s)$. These $a$-points $\delta_a$ are clustered around the critical line $1/2+i\sR$ which happens to be a Julia line for the essential singularity of $\zeta$ at infinity. We observe a remarkable average behaviour for the sequence of values $\zeta(\delta_a)$.    
\end{abstract}

{\small \noindent {\sc Keywords:} Riemann zeta-function, value-distribution, critical line, Julia line}
\bigskip

\section{Motivation and Statement of the Main Results}

In 1879, \'Emile Picard \cite{picard} proved that if an analytic function $f$ has an essential singularity at a point $\omega$, then $f(s)$ takes all possible complex values with at most a single exception (infinitely often) in every neighbourhood of $\omega$. In 1924, Gaston Julia \cite{julia} achieved his celebrated refinement of Picard's great theorem, namely that one can even add a restriction on the angle at $\omega$ to lie in an arbitrarily small cone (see also Burckel \cite[\S 12.4]{burckel}). These results can easily be reformulated for meromorphic functions. Moreover, they can be discussed with respect to functions defined by a Dirichlet series in some half-plane. The latter topic had been investigated to some extent by Szolem Mandelbrojt (in, for example, \cite{mandelbrojt}). In his and later research work (e.g., by Chuji Tanaka \cite{tanaka}) general Dirichlet series $\sum_n a_n\exp(-s\lambda_n)$ are considered and, if certain conditions on the divergent sequence of real exponents $\lambda_n$ are fulfilled, the values are taken around a (horizontal) so-called Julia line. 

In number theory, however, ordinary Dirichlet series appear naturally. The most simple example is the famous Riemann zeta-function, for $\Re s>1$ defined by 
$$
\zeta(s)=\sum_{n\geq 1}n^{-s},
$$
and by analytic continuation elsewhere except for a simple pole at $s=1$. It is well-known that the value-distribution of $\zeta$ is intimately related with the distribution of prime numbers and the yet unsolved Riemann hypothesis on the distribution of zeros would imply the best possible possible error term in the prime number theorem; for this and more we refer to Titchmarsh \cite{tit}. In 1975, Sergey Voronin \cite{voronin} proved with his so-called universality theorem that, roughly speaking, every admissible target function $f$ can be approximated uniformly by appropriate shifts of $\zeta$; more precisely, given a non-vanishing analytic $f$ defined on a closed disk of radius $r\in(0,1/4)$ and $\epsilon>0$, there exists a real number $\tau>0$ such that  
\begin{align*}
\max_{\vert s\vert\leq r}\vert \zeta\left(s+{\textstyle{\frac{3}{4}}}+i\tau\right)-f(s)\vert<\epsilon;
\end{align*}
furthermore, the set of shifts $\tau$ satisfying the above inequality has positive lower density. It is straightforward to deduce that {\it every vertical line $\sigma+i\R$ with $\sigma\in [1/2,1]$ is a Julia line for the essential singularity of $\zeta$ at infinity}. More precisely, {\it for every $\epsilon>0$, every $T>0$ and every complex number $a$ with at most one exception, there exists some $z$ satisfying $\zeta(z)=a$ with $\vert \Re(z)-\sigma\vert<\epsilon$ and $\Im(z)>T$.} Indeed, this may be deduced from the universality property applied to the function constant $a$ (in case $a\neq 0$) and an application of Rouch\'e's theorem (as explained in \cite[\S 8.1 \& 10.1]{steudi}). It should be mentioned that the latter statement about the Julia lines $\sigma+i\R$ for $\sigma>1/2$ also follows from weaker results due to Harald Bohr \& B\"orge Jessen \cite{bohr} from 1932 who showed that $\zeta(\sigma+i\R)$ is dense in $\C$ for $\sigma\in(1/2,1]$. For $\Re s>1$ the defining Dirichlet series converges absolutely and, consequently, the values of $\zeta$ on such vertical lines $\sigma+i\R$ cannot be dense; for $\sigma<1/2$ there is also no such denseness provided that the Riemann hypothesis is true (as shown by Garunk\v{s}tis \& Steuding \cite{garst}), and the remaining case of the critical line is a long-standing folklore conjecture: 
\begin{equation}\label{folk}
\overline{\zeta(1/2+i\R)}=\C\quad ? 
\end{equation}
With regard to Voronin's universality theorem it is known that the range of universality, the open strip $1/2<\Re s<1$, cannot be extended.  

Given an arbitrary complex number $a$, it is an immediate implication of the Julia line $1/2+i\R$ (or universality) that there exists a sequence of complex numbers with real parts converging to $1/2$ on which $\zeta$ assumes the value $a$. For a meromorphic function $f$, the solutions of the equation
$$
f(s)=a
$$
are called $a$-points of $f$, and in the case of the zeta-function investigations of the distribution of $a$-points date back to Edmund Landau's talk \cite{landau} at the International Congress of Mathematicians 1912. It is interesting that a property similar to the aforementioned implication of the Julia line holds for the $a$-points of a different function. Let
\begin{equation}\label{Delta}
\Delta(s):=2^s\pi^{s-1}\Gamma(1-s)\sin{\textstyle{\pi s\over 2}}.
\end{equation}
In this note we shall prove that {\it for every complex number $a\neq 0$ the mean of the values $\zeta(\delta_a)$ on the sequence of $a$-points $\delta_a$ of the function $\Delta$ exists and equals $a+1$}; the case $a=0$ is related to the trivial $\zeta$-zeros at $s=-2n$ for $n\in\N$. Hence, these averages of these $\zeta$-values attain all but one possible complex limit. This indicates an interesting link between the distribution of $a+1$-points of the zeta-functions and $a$-points of $\Delta$. 

The main reason for this remarkable behaviour is the well-known functional equation
\begin{equation}\label{feq}
\zeta(s)=\Delta(s)\zeta(1-s),
\end{equation}
where $\Delta$ is given in (\ref{Delta}). It follows that $\Delta(s)\Delta(1-s)=1$, hence $\Delta({1\over 2}+it)$ lies on the unit circle for real $t$. This special case $a=\exp(2i\phi)$ with real $\phi$ has been studied by Justas Kalpokas and the first author \cite{Kalp} as well as by Kalpokas, Maxim Korolev and the first author \cite{korolev}, respectively. In this situation the $\exp(2i\phi)$-points correspond to intersections of the curve $t\mapsto \zeta({1\over 2}+it)$ with straight lines $\exp(i\phi)\R$ through the origin and related, so-called $\Omega$-results have been derived for these sequences. 

In this note we shall prove for the general case

\begin{Theorem}\label{Th1}
Given a complex number $a\neq 0$, the number $N_a(T)$ of $a$-points $\delta_a=\beta_a+i\gamma_a$ satisfying $\beta_a>-1, 0<\gamma_a<T$ is asymptotically equal to
$$
N_a(T)={T\over 2\pi}\log{T\over 2\pi e}+O_a(\log T),
$$
where the implicit constant in the error term depends on $a$.
\end{Theorem}

\noindent The proof of this theorem will also show that for $a\neq 0$ the $a$-points of $\Delta$ are clustered around the critical line or, in other words, with increasing imaginary part $\gamma_a$ the real part $\beta_a$ is tending to $1/2$. Hence, {\it the critical line $1/2+i\R$ is the unique vertical Julia line for $\Delta$.} There are further $a$-points of $\Delta$ in the left half-plane, close to zeros of $\Delta$ (and they can be located by another application of Rouch\'e's theorem); the condition $\beta_a>-1$ excludes them with at most finitely many exceptions. Notice that $\Delta(s)$ is regular except for simple poles at the positive odd integers $s=2n+1, n\in\N_0$; moreover, $\Delta(s)$ vanishes exactly for the non-positive even integers $s=-2n, n\in\N_0$ (as follows from the well-known properties of Euler's $\Gamma$ and the sine function). Both, $0$ and $\infty$ are thus deficient values for $\Delta$ in the language of value-distribution theory (each of which having defect $\delta=1$). It appears that the distribution of values of both, $\Delta$ and $\zeta$ in the left half-plane is pretty similar (except for the value $0$). In this context the formula in Theorem \ref{Th1} should be compared with the (in principle) identical counterpart for $\zeta$ (obtained by Landau in \cite{andau}).  

For the values of the zeta-function on the $a$-points of $\Delta$ we have

\begin{Theorem}\label{Th2}
For every complex number $a\neq 0$, and every $\eta\in(-\epsilon,1)$, 
\begin{eqnarray*}
\sum_{0<\gamma_a<T\atop \beta_a>-1}\zeta(\eta+\delta_a)&=&{T\over 2\pi}\log{T\over 2\pi e}+{a\over 1-\eta}\left({T\over 2\pi}\right)^{1-\eta}\log{T\over 2\pi} \\
&&-{a\over (1-\eta)^2}\left({T\over 2\pi}\right)^{1-\eta}+O_a(T^{1/2+\epsilon}),
\end{eqnarray*}
where again the implicit constant in the error term depends on $a$ and $\epsilon>0$ is arbitrary.
\end{Theorem}

\noindent Theorem \ref{Th2} implies that $\zeta$ has an essential singularity at infinity (since the values $a$ may be chosen with an arbitrarily big and arbitrarily small absolute value, too big for a removable singularity or a pole). The most interesting case appears for the critical line ($\eta=0$). In this case it follows that 
$$
\lim_{T\to\infty}N_a(T)^{-1}\sum_{0<\gamma_a<T\atop \beta_a>-1}\zeta(\delta_a)=a+1.
$$
The second main term reflects the order of growth of $\zeta(\sigma+it)$ with respect to $\sigma$ (see (\ref{ioi}) and (\ref{oio}) below). In view of $\zeta(\sigma+it)\to 1$ as $\sigma\to +\infty$ the main term tends to $N_a(T)$ for increasing $\eta$. It is certainly possible to enlarge the range for $\eta$, however, we do not consider this of interest.    

\section{Proof of Theorem \ref{Th1}}

By Stirling's formula, 
\begin{equation}\label{pow}
\Delta(\sigma+it)=\left({t\over 2\pi}\right)^{1/2-(\sigma+it)}\exp\left(i(t+{\textstyle{\pi\over 4}})\right)\left(1+O(t^{-1})\right)
\end{equation}
for $t\geq 1$ and $\sigma\geq -2$ (see, for example, \cite[A.2]{neva}). Hence, writing
$$
a=\Delta(\delta_a)=\vert a\vert \exp(i\phi)
$$
for an $a$-point $\delta_a=\beta_a+i\gamma_a$ of $\Delta$, it follows that 
\begin{eqnarray}\label{amod}
\vert a\vert &=& \left({\gamma_a\over 2\pi}\right)^{1/2-\beta_a} \left(1+O(\gamma_a^{-1})\right),\\
\phi &\equiv& \gamma_a\log {2\pi e\over \gamma_a}+{\textstyle{\pi\over 4}}+O(\gamma_a^{-1})\ \bmod\, 2\pi.\nonumber
\end{eqnarray}
This shows that 
\begin{equation}\label{lim}
\beta_a\ \to\ 1/2\qquad \mbox{as}\quad \gamma_a\to\infty. 
\end{equation}
Hence, there exists a positive real number $t_a>1$, depending only on $a$, such that all $a$-points $\delta_a=\beta_a+i\gamma_a$ have a real part $\beta_a\in(-1,2)$ whenever $\gamma_a>t_a$. 

For the proof of Theorem \ref{Th1} we apply the argument principle to the function $\Delta(s)-a$ and integrate counterclockwise over the rectangular contour $\mathcal{C}$ determined by the vertices $-1+it_a$, $2+it_a$, $2+iT$ and $-1+iT$. Hence,
\begin{equation}\label{arg_prin}
2\pi i \left( N_a(T) - N_a(t_a) \right)=\int_{\mathcal{C}} {\Delta'(s)\over \Delta(s)-a} \d s;
\end{equation}
observe that the contour avoids the poles and zeros of $\Delta$ (on the real line). Hence, we may rewrite the integrand as
\begin{equation}\label{bibo}
{\Delta'(s)\over \Delta(s)-a}={\Delta'\over \Delta}(s)\cdot \frac{1}{1 - a/\Delta(s)}.
\end{equation}
Taking into account (\ref{pow}) in combination with a relative of Stirling's formula,
\begin{equation} \label{log_der}
{\Delta' \over \Delta}(\sigma+it) = -\log{t\over 2\pi} + O(t^{-1}),
\end{equation}
which is also valid for $t\geq 1$, we obtain, for $\sigma>1/2$,
\begin{equation}\label{one}
{\Delta'(s)\over \Delta(s)-a}\ll_a t^{1/2-\sigma}\log(t+1),
\end{equation}
where we write $s=\sigma+it$, as is usual in analytic number theory. In a similar fashion we find for $\sigma<1/2$ that
\begin{eqnarray}\label{two}
{\Delta'(s) \over \Delta(s)-a}&=&{\Delta'\over \Delta}(s)\cdot \left(1+\sum_{j\geq 1}\left({a\over \Delta(s)}\right)^j\right)\\
&=& -\log{t\over 2\pi} + O(t^{-1})+O_a\left(t^{\sigma-1/2}\log(t+1) \right)\nonumber;
\end{eqnarray}
expanding the second factor here into a geometric series is justified by (\ref{pow}). 

Using estimate (\ref{one}), the contribution of the integral over the line segments where the integration variable has real part larger than $1/2$ yields
$$ 
\left\{\int_{1/2+it_a}^{2+it_a}+\int_{2+it_a}^{2+iT}+\int_{2+iT}^{1/2+iT}\right\}{\Delta'(s)\over \Delta(s)-a} \d s \ll_a \log T.
$$
And using estimate (\ref{two}), the contribution of the integral over the line segments where the integration variable has real part less than $1/2$ leads to
\begin{eqnarray*}
\lefteqn{\left\{\int_{1/2+iT}^{-1+iT}+\int_{-1+iT}^{-1+it_a}+\int_{-1+it_a}^{1/2+it_a}\right\}{\Delta'(s)\over \Delta(s)-a} \d s}\\
&=& -i\int_{t_a}^T\left(-\log{t\over 2\pi}+O(t^{-1})\right)\d t+O_a(\log T).
\end{eqnarray*}
Therefore, we end up with
$$
\int_{\mathcal{C}} {\Delta'(s)\over \Delta(s)-a} \d s=i T\log{T\over 2\pi e}+O_a(\log T).
$$
Substituting this and the trivial bound $N_a(t_0) = O_a (1)$ into \eqref{arg_prin} finishes the proof of Theorem \ref{Th1}. 

\section{Proof of Theorem \ref{Th2}}

In view of Theorem \ref{Th1} the sequence of imaginary parts of the $a$-points $\delta_a=\beta_a+i\gamma_a$ cannot be too dense. As a matter of fact, for any given $T_0$ there exists $T\in[T_0,T_0+1)$ such that
\begin{equation}\label{condi}
\min_{\delta_a}\vert T-\gamma_a\vert\gg_a 1/\log T,
\end{equation}
where the minimum is being taken over all $\delta_a$. As indicated in the proof of Theorem \ref{Th1} there exists a real number $t_a>1$ such that $\beta_a\in(-1,2)$ for $\gamma_a>t_a$. Hence, the $a$-points in the sum of Theorem \ref{Th2} lie in the interior of the rectangle with vertices $2+it_a,2+iT,-1+iT,-1+it_a$, where $T$ is supposed to satisfy (\ref{condi}). For technical reasons, however, we integrate over the counterclockwise oriented contour ${\mathcal C}$ with vertices $1+\eta+\epsilon+it_a, 1+\eta+\epsilon+iT, -\eta-\epsilon+iT, -\eta-\epsilon+it_a$, where we may assume that $\epsilon<1-\eta$ (so that $1+\eta+\epsilon<2$). In view of (\ref{lim}) the difference is bounded. Hence, 
\begin{eqnarray}\label{i1234}
\lefteqn{\sum_{0<\gamma_a<T\atop \beta_a>-1}\zeta(\eta+\delta_a)}\nonumber \\
&=&\frac{1}{2\pi i}\left\{\int_{1+\eta+\epsilon+it_a}^{1+\eta+\epsilon+iT}+\int_{1+\eta+\epsilon+iT}^{-\eta-\epsilon+iT}+\right.\nonumber\\
&&\qquad \left.+\int_{-\eta-\epsilon+iT}^{-\eta-\epsilon+it_a}+\int_{-\eta-\epsilon+it_a}^{1+\eta+\epsilon+it_a}\right\} \frac{\Delta'(s)}{\Delta(s)-a}\zeta(\eta+s)\d s+O_a(1)\nonumber \\
&=&\sum_{1\leq j\leq 4}\mathcal{I}_j +O_a(1),
\end{eqnarray}
say. For the calculation of these integrals we shall use the standard estimates 
\begin{equation}\label{ioi}
\zeta(\sigma+it)\ll 1+t^{(1-\sigma)/2+\epsilon}\qquad\mbox{for}\quad \sigma\in[0,2],\ t\geq 1,
\end{equation}
and 
\begin{equation}\label{oio}
\zeta(\sigma+it)\ll t^{1/2-\sigma+\epsilon}\qquad\mbox{for}\quad \sigma\in [-1,0],\ t\geq 1
\end{equation}
(see \cite[\S 5.1]{tit}).

We first consider the vertical integrals in (\ref{i1234}). The range of integration of $\mathcal{I}_1$ lies in the half-plane $\sigma>1/2$. Taking into account (\ref{one}) in combination with (\ref{ioi}), we find
\begin{eqnarray}\label{i1}
{\mathcal I}_1 &\ll_a & \int_{t_a}^{T} t^{-1/2-\eta-\epsilon}\log t \d t\ll_a T^{1/2-\eta-\epsilon}\log T.
\end{eqnarray}
For the integral $\mathcal{I}_3$ lying in the half-plane $\sigma<1/2$ we use (\ref{two}) and $\Delta(1-s)\Delta(s)=1$ and get
\begin{eqnarray}\label{i3}
{\mathcal I}_3 &=& -{1\over 2\pi i}\int_{-\eta-\epsilon+it_a}^{-\eta-\epsilon+iT} {\Delta'\over \Delta}(s)\left(1+{a\over \Delta(s)}+\sum_{j\geq 2}\left({a\over \Delta(s)}\right)^j\right)\zeta(\eta+s)\d s\nonumber\\
&=&{1\over 2\pi i}\int_{1+\eta+\epsilon-it_a}^{1+\eta+\epsilon-iT} {\Delta'\over \Delta}(1-s)\left(1+a\Delta(s)+\sum_{j\geq 2}(a\Delta(s))^j\right)\zeta(1+\eta-s)\d s\\
&=&\sum_{1\leq \ell\leq 3}\mathcal{J}_\ell,\nonumber
\end{eqnarray}
say. In general, we observe that
\begin{eqnarray*}
{1\over 2\pi i}\int_{\sigma-it_a}^{\sigma-iT}f_a(s)\d s&=&-{1\over 2\pi}\int_{t_a}^Tf_a(\sigma-it)\d t\\
&=&-{1\over 2\pi}{i\over i}\int_{t_a}^T\overline{f_{\overline{a}}(\sigma+it)}\d t =\overline{-{1\over 2\pi i}\int_{\sigma+it_a}^{\sigma+iT}f_{\overline{a}}(s)\d s}.
\end{eqnarray*}
Thus, substituting the functional equation (\ref{feq}) yields
\begin{eqnarray}\label{grobi}
\overline{\mathcal J}_1 &=& -{1\over 2\pi i}\int_{1+\eta+\epsilon+it_a}^{1+\eta+\epsilon+iT}{\Delta'\over \Delta}(1-s)\Delta(1+\eta-s)\zeta(s-\eta)\d s\nonumber\\
&=&\int_{t_a}^{T}\left(\log\frac{t}{2\pi}+O(t^{-1})\right)\,\d\left(\frac{1}{2\pi i}\int_{1+\eta+\epsilon+i}^{1+\eta+\epsilon+it}\Delta(1+\eta-s)\zeta(s-\eta)\d s\right).
\end{eqnarray}
For the inner integral we use Gonek's lemma: {\it Suppose that $\sum_{n=1}^{\infty}a(n)n^{-s}$ converges for $\sigma>1$ where $a(n)\ll n^\epsilon$ for any $\epsilon >0$. Then we have, uniformly for $1<\alpha \leq 2$,}
\begin{align*}
&
\frac{1}{2\pi i}\int_{\alpha+i}^{\alpha+iT}\left(\frac{m}{2\pi}\right)^s\Gamma(s)\exp\left(\delta \frac{\pi i s}{2}\right)\sum_{n=1}^{\infty}\frac{a(n)}{n^s}\d s\\
&\quad{}=
\left\{\begin{array}{ll}
\sum_{n\leq \frac{Tm}{2\pi}}a(n)\exp(-2\pi i\frac{n}{m})+O\left(m^\alpha T^{\alpha-\frac12+\epsilon}\right)&\mbox{ if } \quad \delta = -1,\\
O(m^a) & \mbox{ if }\quad \delta = +1.
\end{array}\right.
\end{align*}
A proof of this lemma can be found in \cite{gonek} (and there is a pre-version in \cite[\S 7.4]{tit}). Since
$$
\Delta(1-s)=(2\pi)^{-s}\Gamma(s)\left(\exp\left({\pi is\over 2}\right)+\exp\left(-{\pi is\over 2}\right)\right)
$$
by the definition (\ref{Delta}) of $\Delta$, this gives via the substitution $z=s-\eta$ here
$$
\frac{1}{2\pi i}\int_{1+\epsilon+i}^{1+\epsilon+it}\Delta(1-z)\zeta(z)\d z={t\over 2\pi}+O(t^{1/2+\epsilon}).
$$
Plugging this into (\ref{grobi}) leads to
\begin{equation}\label{j1}
\mathcal{J}_1=\frac{T}{2\pi}\log\frac{T}{2\pi e}+O_a\left(T^{1/2+\epsilon}\right).
\end{equation}
We observe that there is no impact of the conjugation in (\ref{grobi}) here and the same will turn out for the further integrals ${\mathcal J}_\ell$ in the sequel.

For $\mathcal{J}_2$ we get   
\begin{eqnarray*}
{\mathcal J}_2 &=& -{a\over 2\pi i}\int_{1+\eta+\epsilon+it_a}^{1+\eta+\epsilon+iT}{\Delta'\over \Delta}(1-s){\Delta(1+\eta-s)\over \Delta(1-s)}\zeta(s-\eta)\d s\\
&=&{a\over 2\pi}\int_{t_a}^{T}\left(\log\frac{t}{2\pi}+O(t^{-1})\right)\left(\left({t\over 2\pi}\right)^{-\eta}+O(t^{-1-\eta})\right)\sum_{n\geq 1}n^{-(1+\epsilon+it)}\d t,
\end{eqnarray*}
where we have used (\ref{pow}) to write
$$
{\Delta(1+\eta-s)\over \Delta(1-s)}=\left({t\over 2\pi}\right)^{-\eta}\left(1+O(t^{-1})\right).
$$
Extracting the constant term (for $n=1$) from the absolutely convergent series yields
\begin{eqnarray}\label{kruemel}
{\mathcal J}_2 &=& {a\over 2\pi}\int_1^T\left({t\over 2\pi}\right)^{-\eta}\log{t\over 2\pi}\d t+O_a(1)+\nonumber\\
&&+O\left(\sum_{n\geq 2}n^{-1-\epsilon}\int_{t_a}^T\left({t\over 2\pi}\right)^{-\eta}\log {t\over 2\pi}\exp(-it\log n)\d t\right).
\end{eqnarray}
For the exponential integral in the latter error term we may use the following lemma: {\it Given real functions $F$ and $G$ on $[a,b]$ such that $G(t)/F'(t)$ is monotonic and $F'(t)/G(t)\geq M>0$ or $F'(t)/G(t)\leq -M<0$, then}
$$
\int_a^bG(t)\exp(iF(t))\d t\ll {1\over M}.
$$
This is essentially a classical lemma from \cite[\S 4.3]{tit}. Taking into account (\ref{kruemel}) this leads in our situation to 
\begin{equation}\label{j2}
{\mathcal J}_2={a \over 1-\eta}\left({T\over 2\pi}\right)^{1-\eta}\log{T\over 2\pi}-{a\over (1-\eta)^2}\left({T\over 2\pi}\right)^{1-\eta}+O_a(1).  
\end{equation}

For the vertical integrals it remains to estimate the integral over the tail of the geometric series. For this purpose we apply (\ref{log_der}) and (\ref{oio}) and get 
$$
\mathcal{J}_3 \ll_a \int_{t_a}^{T} \log t\sum_{j\geq 2}t^{-j(1/2+\eta+\epsilon)}t^{+1/2+\epsilon}\d t \ll_a T^{1/2+\epsilon}.
$$
Substituting this, \eqref{j1} and \eqref{j2} into \eqref{i3}, in combination with \eqref{i1} we conclude that the vertical integrals contribute
\begin{eqnarray}\label{verti}
\mathcal{I}_1+\mathcal{I}_3&=&{T\over 2\pi}\log{T\over 2\pi e}+{a\over \eta-1}\left({T\over 2\pi}\right)^{1-\eta}\log{T\over 2\pi} \\
&&-{a\over (1-\eta)^2}\left({T\over 2\pi}\right)^{1-\eta}+O_a(T^{1/2+\epsilon}),\nonumber  
\end{eqnarray}
which is already the main term. It remains to consider the horizontal integrals in (\ref{i1234}).

Of course, $\mathcal{I}_4\ll_a 1$ (since it is independent of $T$). For $\mathcal{I}_2$, however, we need a little help from Jacques Hadamard's theory of functions of finite order (see \cite[\S 5.3]{palka}). Our approach is pretty similar to the case of the zeta-function (in \cite[\S 9.6]{tit}).

As mentioned above, $\Delta(s)$ is analytic except for simple poles at the positive odd integers. Thus, 
$$
f(s):=(\Delta(s)-a)\cdot \Gamma\left({1-s\over 2}\right)^{-1} 
$$
defines an entire function. It follows from Stirling's formula that $f$ is entire of order one. Hence, Hadamard's factorization theorem yields
$$
f(s)=\exp(A+Bs)\prod_{\delta_a}\left(1-{s\over \delta_a}\right)\exp\left({s\over \delta_a}\right),
$$
where $A$ and $B$ are certain complex constants and the product is taken over {\it all} zeros $\delta_a$ of $f(s)$. Hence, taking the logarithmic derivative, we get
$$
{f'\over f}(s)=B+\sum_{\delta_a}\left(\frac1{s-\delta_a}+\frac1{\delta_a}\right)
$$
In view of 
$$
{f'\over f}(s)={\Delta'(s)\over \Delta(s)-a}+{1\over 2}{\Gamma'\over \Gamma}\left({1-s\over 2}\right)
$$
and 
$$
{\Gamma'\over \Gamma}\left({1-s\over 2}\right)\ll \log t
$$
(as follows from Stirling's formula), we arrive at 
$$
{\Delta'(s)\over \Delta(s)-a}=\sum_{\delta_a}\left(\frac1{s-\delta_a}+\frac1{\delta_a}\right)+O(\log t).
$$
Using this with $s=2+it$ and taking into account that 
$$
{\Delta'(2+it)\over \Delta(2+it)-a}\ll_a 1,
$$
as follows from (\ref{one}), leads to
$$
{\Delta'(s)\over \Delta(s)-a}=\sum_{\delta_a}\left(\frac1{s-\delta_a}-{1\over 2+it-\delta_a}\right)+O_a(\log t).
$$
By Theorem \ref{Th1}, 
$$
\sum_{\vert t-\gamma_a\vert\leq 1}\frac1{2+it-\delta_a}\ll_a \sum_{\vert t-\gamma_a\vert\leq 1}1=N_a(t+1)-N_a(t-1)\ll_a \log t.
$$
Moreover, for any positive integer $n$,
$$
\sum_{t+n<\gamma_a\leq t+n+1}\left(\frac1{s-\delta_a}-{1\over 2+it-\delta_a}\right)\ll_a \sum_{t+n<\gamma_a\leq t+n+1} {1\over n^2}\ll_a {\log (t+n)\over n^2}
$$
as a short computation shows. Since $\sum_{n\geq 1}\log(t+n)/n^2\ll \log t$ it follows that
$$
\sum_{\gamma_a>t+1}\left(\frac1{s-\delta_a}-{1\over 2+it-\delta_a}\right)\ll_a \log t
$$
and by a similar reasoning we can estimate the sum over the $a$-points $\delta_a$ satisfying $\gamma_a<t-1$ by the same bound. Hence, the contribution of the $a$-points distant from $s$ are negligible. This leads to
$$
{\Delta'(s)\over \Delta(s)-a}=\sum_{\vert t-\gamma_a\vert\leq 1}\frac1{s-\delta_a}+O_a(\log t),
$$
valid for $\sigma\in[-1,2]$. 

Hence, the remaining integral $\mathcal{I}_2$ in (\ref{i1234}) can be bounded by
$$
{\mathcal I}_2= {1\over 2\pi i}\int_{+1\eta+\epsilon+iT}^{-\eta-\epsilon+iT}\left(\sum_{\vert T-\gamma_a\vert\leq 1}\frac1{s-\delta_a}+O_a(\log T)\right)\zeta(\eta+s)\d s. 
$$
In view of (\ref{condi}) we have $1/\vert s-\delta_a\vert\ll_a \log T$. Thus, using Theorem \ref{Th1}, (\ref{oio}) and (\ref{ioi}) yields  
\begin{eqnarray*}
{\mathcal I}_2&\ll_a &(\log T)^2 \left\{\int_{-\epsilon+iT}^{+iT}+\int_{+iT}^{1+iT}+\int_{1+iT}^{1+2\eta+\epsilon+iT}\right\}\vert \zeta(s)\vert\d s\\
&\ll_a & (\log T)^2\left\{T^{1/2+\epsilon}+T^{1/2+\epsilon}+1\right\}\ll_a T^{1/2+\epsilon}, 
\end{eqnarray*}
where $\epsilon$ at different places may take different values. This is the contribution of the horizontal integrals. In combination with (\ref{verti}) we thus arrive via (\ref{i1234}) at the asymptotic formula of the theorem. In view of (\ref{ioi}) and (\ref{oio}) we may replace the chosen $T$ (with respect to (\ref{condi}) by a general $T$ at the expense of an error of order $T^{1/2+\epsilon}$ (as follows from estimate (\ref{ioi})). This finally proves Theorem \ref{Th2}.  

\section{Concluding Remarks}

Our results rely heavily on the functional equation (\ref{feq}) and its implication $\Delta(s)\Delta(1-s)=1$. These identities imply a certain point symmetry for both $\zeta$ and $\Delta$ with respect to $s=1/2$ which is the abscissa of the vertical Julia line for $\Delta$. Concerning the location of the $a$-points of $\Delta$ let $a\in\C$ be non-zero. All complex numbers of absolute value $\vert a\vert$ form a circle of radius $\vert a\vert$ centered at the origin. According to (\ref{amod}) the $a$-points of $\Delta$ lie close to the curves $C_{\vert a\vert}$ defined by 
$$
\vert a\vert = \left({t\over 2\pi}\right)^{1/2-\sigma}
$$
or
$$
\sigma=1/2-{\log\vert a\vert\over \log (t/(2\pi)},
$$
where we may assume that $t$ is sufficiently large. With increasing real part these $\delta_a=\beta_a+i\gamma_a$ lie closer and closer to $C_{\vert a\vert}$. It follows that circles of radius $\vert a\vert<1$ correspond to curves $C_{\vert a\vert}$ to the right of the Julia line $1/2+i\R$ and circles of radius $\vert a\vert>1$ correspond to curves $C_{\vert a\vert}$ to the left. Asymptotically all these curves approach the line $1/2+i\R$ which is the exact curve corresponding to the unit circle or the points $a$ with $\vert a\vert=1$.  

It follows from the two theorems that 
$$
\zeta(1/2+i\gamma_a)-\zeta(\delta_a)=\int_{\beta_a}^{1/2}\zeta'(\sigma+i\gamma_a)\d\sigma,
$$
hence
\begin{eqnarray*}
\sum_{0<\gamma_a<T\atop \beta_a>-1}\zeta(1/2+i\gamma_a)&=&(a+1)N_a(T)+O_a(T^{1/2+\epsilon})+\\
&&+O\left(\sum_{0<\gamma_a<T\atop \beta_a>-1}\left\vert \beta_a-1/2\right\vert\max_{1/2\leq \sigma\leq \beta_a\atop \mbox{\small or}\,\beta_a\leq \sigma\leq 1/2}\vert \zeta'(\sigma+i\gamma_a)\vert\right).
\end{eqnarray*}
One could imagine that also the values $\zeta(1/2+i\gamma_a)$ equal $a+1$ on average but for a proof too little is known about the quantities $\zeta'(\sigma+i\gamma_a)$. The counterpart to the existence of the mean of the values $\zeta(\delta_a)$ with respect to the point symmetry at $s=1/2$ given by (\ref{feq}) is that $\zeta(1-\delta_a)$ is on average equal to $1+1/a$. And for results with respect to $a$-points of $\Delta$ in the lower half-plane one may use the functional equation (\ref{feq}) in combination with the reflection principle $\zeta(\overline{s})=\overline{\zeta(s)}$. Unfortunately, our results do not shed new light on the question proposed by (\ref{folk}). 

It is also straightforward to prove similar results for other functions defined by a Dirichlet series in some half-plane and satisyfing a Riemann-type functional equation (for example Dirichlet $L$-functions to residue class characters).

\bigskip

\noindent
J\"orn Steuding\\
Department of Mathematics, W\"urzburg University\\ 
Emil Fischer-Str. 40, 97\,074 W\"urzburg, Germany\\
steuding@mathematik.uni-wuerzburg.de
\medskip

\noindent 
Ade Irma Suriajaya\\
Faculty of Mathematics, Kyushu University \\
744 Motooka, Nishi-ku, Fukuoka 819-0395, Japan\\
adeirmasuriajaya@math.kyushu-u.ac.jp

\bigskip

\end{document}